\documentclass{amsart}
\usepackage{latexsym}
\usepackage{amssymb}
\usepackage{amsmath}
\def\gg{\mathfrak{G}}
\def\gk{\mathfrak{K}}
\def\gp{\mathfrak{P}}
\def\gt{\mathfrak{T}}
\newtheorem{df}{Definition}[section]

\newtheorem{theo}[df]{Theorem}
\newtheorem{remark}[df]{Remark}

\newtheorem{lem}[df]{Lemma}

\begin{document}

\title[The First Eigenvalue of the Dirac Operator...]{The First Eigenvalue of the Dirac Operator on Compact Spin Symmetric Spaces}%
\author{Jean-Louis Milhorat}
\address{Laboratoire Jean Leray\\
 UMR CNRS 6629\\
 D\'epartement de Math\'ematiques\\
 Universit\'e de Nantes\\
 2, rue de la Houssini\`ere\\
 BP 92208\\
  F-44322 NANTES CEDEX 03}
\email{milhorat@math.univ-nantes.fr}


\date{}
\begin{abstract}
We give a formula for the first eigenvalue of the Dirac operator
acting on spinor fields of a spin compact irreducible symmetric
space $G/K$.
\end{abstract}

\maketitle
\section{Introduction}
It is well-known that symmetric spaces provide examples where
detailed information on the spectrum of Laplace or Dirac operators
can be obtained. Indeed, for those manifolds, the computation of
the spectrum can be (theoretically) done using group theoretical
methods. However the explicit computation is far from being simple
in general and only few examples are known.  On the other hand,
many results require some information about the first (nonzero)
eigenvalue, so it seems interesting to get this eigenvalue without
computing all the spectrum. In that direction, the aim of this
paper is to prove the following formula for the first eigenvalue
of the Dirac operator:
\begin{theo} Let $G/K$ be a compact,
simply-connected, $n$-dimensional irreducible symmetric space with
$G$ compact and simply-connected, endowed with the metric induced
by the Killing form of $G$ sign-changed.  Assume that $G$ and $K$
have same rank and that $G/K$ has a spin structure. Let
$\beta_{k}$, $k=1,\ldots, p$, be the $K$-dominant weights
occurring in the decomposition into irreducible components of the
spin representation under the action of $K$. Then the square of
the first eigenvalue of the Dirac operator is
\begin{equation}\label{sqvp}
 2\,\min_{1\leq k\leq p}\|\beta_{k}\|^{2}+n/8\,,
\end{equation}
where $\|\cdot\|$ is the norm associated to the scalar product
$<\,,\,>$ induced by the Killing form of $G$ sign-changed.
\end{theo}
\begin{remark}
The proof uses a lemma of R. Parthasarathy in \cite{Par}, which
allows to express (\ref{sqvp}) in the following way. Let
 $T$ be a fixed common maximal torus of $G$ and $K$. Let $\Phi$ be
  the set of non-zero roots of $G$ with respect to $T$.
Let $\delta_{G}$, (resp. $\delta_{K}$) be the half-sum of the
positive roots of $G$, (resp. $K$), with respect to a fixed
lexicographic ordering in $\Phi$. Then the square of the first eigenvalue of the
 Dirac operator is given by
\begin{equation}
 2\,\|\delta_{G}\|^{2}+2\,
 \|\delta_{K}\|^{2}-4\, \max_{w\in
 W}<w\cdot\delta_{G},\delta_{K}>+n/8\,,
\end{equation}
where $W$ is a certain (well-defined) subset of the Weyl group of
$G$.
\end{remark}

\section{The Dirac Operator on a Spin Compact Symmetric Space}

We first review some results about the Dirac operator on a spin
symmetric space, cf.~for instance \cite{CFG89} or \cite{Bar91}. A
detailed survey on the subject may be found, among other topics,
in the reference \cite{BHMM}. Let $G/K$ be a spin compact
symmetric space. We assume that $G/K$ is simply connected, so $G$
may be chosen to be compact and simply connected and $K$ is the
connected subgroup formed by the fixed elements of an involution
$\sigma$ of $G$, cf. \cite{Hel}. This involution induces the
Cartan decomposition of the Lie algebra $\gg$ of $G$ into
$$\gg=\gk\oplus\gp\,,$$ where $\gk$ is the Lie algebra of $K$ and
$\gp$ is the vector space $\{X\in\gg\,;\, \sigma_{*}\cdot X=-X\}$.
This space $\gp$ is canonically identified with the tangent space
to $G/K$ at the point $o$, $o$ being the class of the neutral
element of $G$. We also assume that the symmetric space $G/K$ is
irreducible, so all the $G$-invariant scalar products on $\gp$,
hence all the $G$-invariant Riemannian metrics on $G/K$ are
proportional. We consider the metric induced by the Killing form
of $G$ sign-changed. With this metric, $G/K$ is an Einstein space
with scalar curvature $\mathrm{Scal}=n/2$. The spin condition
implies that the homomorphism $\alpha: K\rightarrow
\mathrm{SO}(\gp)\simeq \mathrm{SO}_{n}$, $k\mapsto
\mathrm{Ad}_{G}(k)_{|\gp}$ lifts to a homomorphism
$\widetilde{\alpha}:K\rightarrow \mathrm{Spin}_{n}$,
cf.~\cite{CG}. Let $\rho:\mathrm{Spin}_{n}\rightarrow
\mathrm{Hom}_{\mathbb{C}}(\Sigma,\Sigma)$ be the spin
representation. The composition $\rho\circ\widetilde{\alpha}$
defines a ``{spin}'' representation of $K$ which is denoted
$\rho_{K}$. The spinor bundle is then isomorphic to the vector
bundle $$\boldsymbol{\Sigma}:=G\times_{\rho_{K}}\Sigma\,.$$ Spinor
fields on $G/K$ are then viewed as $K$-equivariant functions
$G\rightarrow \Sigma$, i.e. functions:
 $$\Psi:G\rightarrow\Sigma\quad \mathrm{s.t.}\quad \forall g\in G\;,\;\forall k\in K\;,\;
 \Psi(gk)=\rho_{K}(k^{-1})\cdot\Psi(g)\,.$$
Let $L_{K}^{2}(G,\Sigma)$ be the Hilbert space of $L^{2}$
$K$-equivariant functions $G\rightarrow \Sigma$. The Dirac
operator $\mathcal{D}$ extends to a self-adjoint operator on
$L_{K}^{2}(G,\Sigma)$. Since it is an elliptic operator, it has a
(real) discrete spectrum. Now
 if the spinor field $\Psi$ is an eigenvector of $\mathcal{D}$
 for the eigenvalue $\lambda$, then the spinor field $\sigma^{*}\cdot\Psi$
 is an eigenvector for the eigenvalue $-\lambda$, hence the
 spectrum of the Dirac operator is symmetric with respect to the
 origin. Thus the spectrum of $\mathcal{D}$ may be deduced from
 the spectrum of its square $\mathcal{D}^{2}$. By the Peter-Weyl theorem,
 the natural unitary representation of $G$ on the Hilbert space
  $L_{K}^{2}(G,\Sigma)$
 decomposes into the Hilbert sum
 $$ \mathop{\oplus}_{\gamma \in \widehat{G}}V_{\gamma}\otimes
\mathrm{Hom}_{K}(V_{\gamma},\Sigma)\,,$$ where $\widehat{G}$ is
the set of equivalence classes of irreducible unitary complex
representations of $G$, $(\rho_{\gamma},V_{\gamma})$ represents an
element $\gamma\in\widehat{G}$ and
$\mathrm{Hom}_{K}(V_{\gamma},\Sigma)$ is the vector space of
$K$-equivariant homomorphisms $V_{\gamma}\rightarrow \Sigma$, i.e.
$$\mathrm{Hom}_{K}(V_{\gamma},\Sigma)=\{ A\in
\mathrm{Hom}(V_{\gamma},\Sigma)\; \mathrm{s.t.}\; \forall k\in
K\,, A\circ \rho_{\gamma}(k)=\rho_{K}(k)\circ A\}\,.$$ The
injection $ V_{\gamma}\otimes
\mathrm{Hom}_{K}(V_{\gamma},\Sigma)\hookrightarrow
L_{K}^{2}(G,\Sigma)$ is given by $$v\otimes A \mapsto
\Big(g\mapsto (A\circ\rho_{\gamma}(g^{-1})\,)\cdot v\Big)\,.$$
Note that $V_{\gamma}\otimes \mathrm{Hom}_{K}(V_{\gamma},\Sigma)$
consists of $\mathcal{C}^{\infty}$ spinor fields to which the
Dirac operator can be applied. The restriction of
$\mathcal{D}^{2}$ to the space $V_{\gamma}\otimes
\mathrm{Hom}_{K}(V_{\gamma},\Sigma)$ is given by the Parthasaraty
formula, \cite{Par}:
\begin{equation}\label{Par}
\mathcal{D}^{2}(v\otimes A)= v\otimes
(A\circ \mathcal{C}_{\gamma})+\frac{\mathrm{Scal}}{8}\,v\otimes
A\,,
\end{equation} where $\mathcal{C}_{\gamma}$ is the Casimir
operator of the representation $(\rho_{\gamma},V_{\gamma})$. Now
since the representation is irreducible, the Casimir operator is a
scalar multiple of identity, $\mathcal{C}_{\gamma}=c_{\gamma}\,
\mathrm{id}$, where the eigenvalue $c_{\gamma}$ only depends of
$\gamma\in\widehat{G}$. Hence if
$\mathrm{Hom}_{K}(V_{\gamma},\Sigma)\neq\{0\}$, $c_{\gamma}+n/16$
belongs to the spectrum of $\mathcal{D}^{2}$. Let $\rho_{K}=\oplus
\rho_{K,k}$ be the decomposition of the spin representation
$K\rightarrow\Sigma$ into irreducible components. Denote by
$\mathrm{m}({\rho_{\gamma}}_{|K},\rho_{K,k})$ the multiplicity of
the irreducible $K$-representation $\rho_{K,k}$ in the
representation $\rho_{\gamma}$ restricted to $K$. Then $$\dim\,
\mathrm{Hom}_{K}(V_{\gamma},\Sigma)=\sum_{k}
\mathrm{m}({\rho_{\gamma}}_{|K},\rho_{K,k})\,.$$ So the spectrum
of the square of the Dirac operator is
\begin{equation}\label{spec}
\mathrm{Spec}(\mathcal{D}^{2})=\{c_{\gamma}+n/16\;;\;
\gamma\in\widehat{G}\;\mathrm{s.t.}\;\exists k\;\mathrm{s.t.}\;
\mathrm{m}({\rho_{\gamma}}_{|K},\rho_{K,k})\neq
0\}\,.\end{equation}

\section{Proof of the result}
 We assume that $G$ and $K$ have same rank. Let $T$ be a fixed
 common maximal torus. Let $\Phi$ be the set of non-zero roots of the group $G$
 with respect to $T$. According to a classical terminology, a root
 $\theta$ is called compact if the corresponding root space is contained
 in $\gk_{\mathbb{C}}$ (that is, $\theta$ is a root of $K$ with respect to $T$) and
 noncompact if the root space is contained in $\gp_{\mathbb{C}}$.
Let $\Phi_{G}^{+}$  be the set of positive roots of $G$,
$\Phi_{K}^{+}$
 be the set of positive roots of $K$,
and $\Phi_{n}^{+}$ be the set of positive noncompact roots with
respect to a fixed lexicographic ordering in $\Phi$. The half-sums
of the positive roots of $G$ and $K$ are respectively denoted
$\delta_{G}$ and $\delta_{K}$ and the half-sum of noncompact
positive roots is denoted by $\delta_{n}$. The Weyl group of $G$
is denoted $W_{G}$. The space of weights is endowed with the
$W_{G}$-invariant scalar product $<\,,\,>$ induced by the Killing
form of $G$ sign-changed.\\ Let
\begin{equation}W:=\{w\in W_{G}\;;\;
w\cdot \Phi_{G}^{+}\supset \Phi_{K}^{+}\}\,.
\end{equation}
By a result of R. Parthasaraty, cf.~lemma~2.2 in \cite{Par}, the
spin representation $\rho_{K}$ of $K$ decomposes into the
irreducible sum
\begin{equation}
\rho_{K}=\bigoplus_{w\in W}\rho_{K,w}\;,
\end{equation}
where $\rho_{K,w}$ has for dominant weight
\begin{equation}
\beta_{w}:= w\cdot \delta_{G}-\delta_{K}\,.
\end{equation}
Now define $w_{0}\in W$ such that
\begin{equation}\label{w0}
 \|\beta_{w_{0}}\|^{2}=\min_{w\in W}\|\beta_{w}\|^{2}\,,
\end{equation}
and
\begin{equation}\label{w1}
 \mbox{if there exists a } w_{1}\neq w_{0}\in W \mbox{ such that }
\|\beta_{w_{1}}\|^{2}=\min_{w\in W}\|\beta_{w}\|^{2}\, ,\mbox {
 then } \beta_{w_{1}}\prec \beta_{w_{0}}\,,
\end{equation}
where $\prec$ is the usual ordering on weights.
\begin{lem} The weight
$$\beta_{w_{0}}^{G}:=w_{0}^{-1}\cdot\beta_{w_{0}}=\delta_{G}-w_{0}^{-1}\cdot\delta_{K}\,,$$
is $G$-dominant.
\end{lem}
\begin{proof} Let $\Pi_{G}=\{\theta_{1},\ldots,\theta_{r}\}\subset \Phi_{G}^{+}$ be the
set of simple roots. It is sufficient to prove that
$2\,\frac{<\beta_{w_{0}}^{G},\theta_{i}>}{<\theta_{i},\theta_{i}>}$
is a non-negative integer for any simple root $\theta_{i}$. Since
$T$ is a maximal common torus of $G$ and $K$, $\beta_{w_{0}}$,
which is an integral weight for $K$ is also an integral weight for
$G$. Now since the Weyl group $W_{G}$ permutes the weights,
$\beta_{w_{0}}^{G}=w_{0}^{-1}\cdot\beta_{w_{0}}$ is also a
integral weight for $G$, hence
$2\,\frac{<\beta_{w_{0}}^{G},\theta_{i}>}{<\theta_{i},\theta_{i}>}$
is an integer for any simple root $\theta_{i}$. So we only have to
prove that this integer is non-negative.\\  Let $\theta_{i}$ be a
simple root. Since
$2\,\frac{<\delta_{G},\theta_{i}>}{<\theta_{i},\theta_{i}>}=1$,
(see for instance \S~10.2 in \cite{Hum}) and since the scalar
product $<\cdot ,\cdot>$ is $W_{G}$-invariant, one gets
\begin{equation}
2\,\frac{<\beta_{w_{0}}^{G},\theta_{i}>}{<\theta_{i},\theta_{i}>}=1-
2\,\frac{<\delta_{K},w_{0}\cdot\theta_{i}>}{<\theta_{i},\theta_{i}>}\,.
\end{equation}
Suppose first that $w_{0}\cdot\theta_{i}\in \Phi_{K}$. If
$w_{0}\cdot\theta_{i}$ is positive then $w_{0}\cdot\theta_{i}$ is
necessarily a $K$-simple root. Indeed let
$\Pi_{K}=\{\theta_{1}',\ldots,\theta_{l}'\}\subset \Phi_{K}^{+}$
be the set of $K$-simple roots. One has
$w_{0}\cdot\theta_{i}=\sum_{j=1}^{l} b_{ij}\, \theta_{j}'$, where
the $b_{ij}$ are non-negative integers. But since $w_{0}\in W$,
there are $l$ positive roots $\alpha_{1},\ldots,\alpha_{l}$ in
$\Phi_{G}^{+}$ such that $w_{0}\cdot \alpha_{j}=\theta_{j}'$,
$j=1,\ldots, l$. So $ \theta_{i}=\sum_{j=1}^{l}
b_{ij}\,\alpha_{j}$. Now each $\alpha_{j}$ is a sum of simple
roots $\sum_{k=1}^{r}a_{jk}\,\theta_{k}$, where the $a_{jk}$ are
non-negative integers. So $\theta_{i}=\sum_{j,k} b_{ij}\,
a_{jk}\,\theta_{k}$. By the linear independence of simple roots,
one gets $\sum_{j} b_{ij}\, a_{jk}=0$ if $k\neq i$, and $\sum_{j}
b_{ij}\, a_{ji}=1$. Hence there exists a $j_{0}$ such that
$b_{ij_{0}}=a_{j_{0}i}=1$, the other coefficients being zero. So
$w_{0}\cdot\theta_{i}=\theta_{j_{0}}'$ is a $K$-simple root. Now
since
$2\,\frac{<\delta_{K},w_{0}\cdot\theta_{i}>}{<\theta_{i},\theta_{i}>}=
2\,\frac{<\delta_{K},w_{0}\cdot\theta_{i}>}
{<w_{0}\cdot\theta_{i},w_{0}\cdot\theta_{i}>}=1$, one gets
$2\,\frac{<\beta_{w_{0}}^{G},\theta_{i}>}{<\theta_{i},\theta_{i}>}=0$,
hence
$2\,\frac{<\beta_{w_{0}}^{G},\theta_{i}>}{<\theta_{i},\theta_{i}>}\geq
0$. Now, the same conclusion holds if $w_{0}\cdot\theta_{i}$ is a
negative root of $K$, since
$2\,\frac{<\delta_{K},w_{0}\cdot\theta_{i}>}{<\theta_{i},\theta_{i}>}=-2\,
\frac{<\delta_{K},-w_{0}\cdot\theta_{i}>}{<\theta_{i},\theta_{i}>}=-1$,
hence
$2\,\frac{<\beta_{w_{0}}^{G},\theta_{i}>}{<\theta_{i},\theta_{i}>}=2$.\\
Suppose now that $w_{0}\cdot\theta_{i}\notin \Phi_{K}$, that is
$w_{0}\cdot\theta_{i}$ is a noncompact root. This implies that
$w_{0}\sigma_{i}$, where $\sigma_{i}$ is the reflection across the
hyperplane $\theta_{i}^{\bot}$, is an element of $W$. Let $
\alpha_{1},\ldots, \alpha_{m}$ be the positive roots in
$\Phi_{G}^{+}$ such that $w_{0}\cdot \alpha_{j}=\alpha'_{j}$,
where the $\alpha'_{j}$, $j=1,\ldots, m$ are the positive roots of
$K$. Since $\sigma_{i}$ permutes the positive roots other than
$\theta_{i}$, (cf. for instance Lemma~B, \S~10.2 in \cite{Hum}),
and since $\theta_{i}$ can not be one of the roots
$\alpha_{1},\ldots,\alpha_{m}$ (otherwise $w_{0}\cdot\theta_{i}\in
\Phi_{K}^{+}$), each root $\sigma_{i}\cdot \alpha_{j}$ is
positive. So $w_{0}\sigma_{i}\in W$ since
$w_{0}\sigma_{i}\cdot(\sigma_{i}\cdot \alpha_{j})=\alpha_{j}'$,
$j=1,\ldots,m$.\\ We now claim that
$2\,\frac{<\beta_{w_{0}}^{G},\theta_{i}>}{<\theta_{i},\theta_{i}>}<0$,
which is equivalent to
$2\,\frac{<\delta_{K},w_{0}\cdot\theta_{i}>}{<\theta_{i},\theta_{i}>}>1$,
is impossible.\\ Suppose that
\begin{equation}
2\,\frac{<\delta_{K},w_{0}\cdot\theta_{i}>}{<\theta_{i},\theta_{i}>}>1\,.
\end{equation}
Since $\delta_{K}$ can be expressed as $\delta_{K}=\sum_{i=1}^{l}
c_{i}\,\theta_{i}'$, where the $c_{i}$ are nonnegative, there
exists a $K$-simple root $\theta_{j}'$ such that
$<\theta_{j}',w_{0}\cdot\theta_{i}>0$, and since
$2\,\frac{<\theta_{j}',w_{0}\cdot\theta_{i}>}{<\theta_{j}',\theta_{j}'>}$
is an integer, this implies that
\begin{equation}\label{cond1}
2\,\frac{<\theta_{j}',w_{0}\cdot\theta_{i}>}{<\theta_{j}',\theta_{j}'>}\geq
1\,.
\end{equation}
So $\theta_{j}'-w_{0}\cdot\theta_{i}$ is a root (cf. for instance
\S~9.4 in \cite{Hum}). Moreover, from the bracket relation
$[\gk,\gp]\subset\gp$, it is a noncompact root. Now $\pm
(\theta_{j}'-w_{0}\cdot\theta_{i})$ is a positive noncompact root,
so by the description of the weights of the spin representation
$\rho_{K}$, (they are of the form: $\delta_{n}-$(a sum of distinct
positive noncompact roots), cf. \S2~in \cite{Par}),
$$(w_{0}\cdot\delta_{G}-\delta_{K})\pm
(\theta_{j}'-w_{0}\cdot\theta_{i}) \mbox{ is a weight of }
\rho_{K}\,.$$  Now, $(w_{0}\cdot\delta_{G}-\delta_{K})+
(\theta_{j}'-w_{0}\cdot\theta_{i})$ can not be a weight of
$\rho_{K}$. Otherwise since
$\sigma_{i}\cdot\delta_{G}=\delta_{G}-\theta_{i}$,
$(w_{0}\sigma_{i}\cdot\delta_{G}-\delta_{K})+ \theta_{j}'$ is a
weight of $\rho_{K}$. But since $w_{0}\sigma_{i}\in W$,
$\mu:=w_{0}\sigma_{i}\cdot\delta_{G}-\delta_{K}$ is a dominant
weight of $\rho_{K}$. So $\mu$ is a dominant weight but not the
highest weight of an irreducible component of $\rho_{K}$. Hence
there exists an irreducible representation of $\rho_{K}$ with
dominant weight $\lambda=w\cdot\delta_{G}-\delta_{K}$, $w\in W$,
whose set of weights $\Pi$ contains $\mu$. Furthermore
$\mu\prec\lambda$. Now since $\mu \in \Pi$$,
\|\mu+\delta_{K}\|^{2}\leq \|\lambda+\delta_{K}\|^{2}$, with
equality only if $\mu=\lambda$, (cf. for instance Lemma~C, \S 13.4
in \cite{Hum}). But $\|\mu+\delta_{K}\|^{2}=\|\delta_{G}\|^{2}=
\|\lambda+\delta_{K}\|^{2}$, so $\mu=\lambda$, contradicting the
fact that $\mu\prec\lambda$.  \\ Thus only
\begin{equation}\mu_{0}:=(w_{0}\cdot\delta_{G}-\delta_{K})-
(\theta_{j}'-w_{0}\cdot\theta_{i})\,, \end{equation} can be a
weight of $\rho_{K}$. Now one has $$
  \begin{array}{rl}
    \|\mu_{0}\|^{2}=
&\|w_{0}\cdot\delta_{G}-\delta_{K}+ w_{0}\cdot\theta_{i}\|^{2}\\
&-2\,<w_{0}\cdot\delta_{G}-\delta_{K}+w_{0}\cdot\theta_{i},\theta_{j}'>
+\|\theta_{j}'\|^{2}\,.
  \end{array}$$
Since $w_{0}\cdot\delta_{G}-\delta_{K}$ is a dominant weight,
$<w_{0}\cdot\delta_{G}-\delta_{K},\theta_{j}'>\geq 0$, and from
(\ref{cond1}), $2\,<w_{0}\cdot\theta_{i},\theta_{j}'>
-\|\theta_{j}'\|^{2}\geq 0$, so $$ \|\mu_{0}\|^{2}\leq
\|(w_{0}\cdot\delta_{G}-\delta_{K})+w_{0}\cdot\theta_{i}\|^{2}\,.$$
Now $$
  \begin{array}{rl}
   \|(w_{0}\cdot\delta_{G}-\delta_{K})+w_{0}\cdot\theta_{i}\|^{2}=
&\|w_{0}\cdot\delta_{G}-\delta_{K}\|^{2}\\
&+2\,<\delta_{G}-w_{0}^{-1}\cdot\delta_{K},\theta_{i}>
+\|\theta_{i}\|^{2}\,.
  \end{array}$$
But, as we supposed
$2\,\frac{<\beta_{w_{0}}^{G},\theta_{i}>}{<\theta_{i},\theta_{i}>}<0$,
one has
$\frac{2\,<\delta_{G}-w_{0}^{-1}\cdot\delta_{K},\theta_{i}>}{\|\theta_{i}\|^{2}}\leq
-1$, so\\ $2\,<\delta_{G}-w_{0}^{-1}\cdot\delta_{K},\theta_{i}>
+\|\theta_{i}\|^{2}\leq 0$, hence
$$\|(w_{0}\cdot\delta_{G}-\delta_{K})+w_{0}\cdot\theta_{i}\|^{2}\leq
\|w_{0}\cdot\delta_{G}-\delta_{K}\|^{2}\,,$$ so
$$\|\mu_{0}\|^{2}\leq \|w_{0}\cdot\delta_{G}-\delta_{K}\|^{2}\,.$$
Now, being a weight of $\rho_{K}$, $\mu_{0}$ is conjugate under
the Weyl group of $K$ to a dominant weight of $\rho_{K}$, say
$w_{1}\cdot\delta_{G}-\delta_{K}$, with $w_{1}\in W$. Note that
$w_{1}\neq w_{0}$, otherwise since $\mu_{0}\prec
w_{1}\cdot\delta_{G}-\delta_{K}$, (cf. Lemma~A, \S~13.2 in
\cite{Hum}), the noncompact root
$\theta_{j}'-w_{0}\cdot\theta_{i}$ should be a linear combination
with integral coefficients of
 compact simple roots. But, by the bracket relation
 $[\gk,\gk]\subset\gk$, that is impossible. Thus, by the definition of $w_{0}$, cf.
(\ref{w0}), $\|w_{0}\cdot\delta_{G}-\delta_{K}\|^{2}\leq
\|w_{1}\cdot\delta_{G}-\delta_{K}\|^{2}=\|\mu_{0}\|^{2}$, so $$
\|\mu_{0}\|^{2}=\|w_{1}\cdot\delta_{G}-\delta_{K}\|^{2}
=\|w_{0}\cdot\delta_{G}-\delta_{K}\|^{2}\,.$$ But by the
condition~(\ref{w1}), the last equality is impossible, otherwise
since $\mu_{0}\prec w_{1}\cdot\delta_{G}-\delta_{K}$ and
$w_{1}\cdot\delta_{G}-\delta_{K}\prec
w_{0}\cdot\delta_{G}-\delta_{K}$, the noncompact root
$\theta_{j}'-w_{0}\cdot\theta_{i}$ should be a linear combination
with integral coefficients of
 compact simple roots. Hence
$2\,\frac{<\beta_{w_{0}}^{G},\theta_{i}>}{<\theta_{i},\theta_{i}>}\geq
0$ also if $w_{0}\cdot\theta_{i}\notin \Phi_{K}$.
\end{proof}

 Now let $(\rho_{0},V_{0})$ be an irreducible representation of
$G$ with dominant weight $\beta_{w_{0}}^{G}$. The fact that
$\beta_{w_{0}}=w_{0}\cdot \beta_{w_{0}}^{G}$ is a weight of
$\rho_{0}$ is an indication that ${\rho_{0}}_{|K}$ may contain the
irreducible representation $\rho_{K,w_{0}}$. This is actually
true:
\begin{lem} With the notations above,
$$\mathrm{m}({\rho_{0}}_{|K},\rho_{K,w_{0}})\geq 1\,.$$
\end{lem}
\begin{proof} Let $v_{0}$ be the maximal vector in $V_{0}$, (it is
unique up to a nonzero scalar multiple). Let $g_{0}\in T$ be a
representative of $w_{0}$. Then $g_{0}\cdot v_{0}$ is a weight
vector for the weight $\beta_{w_{0}}$, since for any $X$ in the
Lie algebra $\gt$ of $T$: $$
\begin{array}{rll}
   X\cdot (g_{0}\cdot v_{0})  &=\frac{d}{dt}\Big((\exp(tX)\, g_{0})\cdot v_{0}\Big)_{|t=0}
   & =\frac{d}{dt}\Big(\Big(g_{0}\,g_{0}^{-1}\exp(tX)\, g_{0}\Big)\cdot v_{0}\Big)_{|t=0}\\
     &= g_{0}\cdot\Big( \Big(\mathrm{Ad}(g_{0}^{-1})\cdot X\Big)\cdot
     v_{0}\Big)&=\beta_{w_{0}}^{G}(w_{0}^{-1}\cdot X)\,(g_{0}\cdot v_{0})\\
     &=(w_{0}\cdot \beta_{w_{0}}^{G})(X)\,(g_{0}\cdot v_{0})&=\beta_{w_{0}}(X)\,
     (g_{0}\cdot v_{0})\,.
\end{array}$$
In order to prove the result, we only have to prove that
$g_{0}\cdot v_{0}$ is a maximal vector (for the action $K$), hence
is killed by root-vectors corresponding to simple roots of $K$. So
let $\theta_{i}'$ be a simple root of $K$ and $E_{i}'$ be a
root-vector corresponding to that simple root. Since $w_{0}\in W$,
there exists a positive root $\alpha_{i}\in \Phi_{G}^{+}$ such
that $w_{0}\cdot\alpha_{i}=\theta_{i}'$. Then $ E_{i}:=
\mathrm{Ad}(g_{0}^{-1})(E_{i}')$ is a root-vector corresponding to
the root $\alpha_{i}$ since for any $X$ in $\gt$ $$
  \begin{array}{rll}
   [X, E_{i}]&= [X,\mathrm{Ad}(g_{0}^{-1})(E_{i}')] &
    =\mathrm{Ad}(g_{0}^{-1})\cdot [\mathrm{Ad}(g_{0})(X),E_{i}'] \\
    &=\mathrm{Ad}(g_{0}^{-1})\cdot[w_{0}\cdot X,E_{i}']
    &=\Big((w_{0}^{-1}\cdot\theta_{i}')(X)\Big)\;\mathrm{Ad}(g_{0}^{-1})\cdot E_{i}'\\
    &=\alpha_{i}(X)\; E_{i}\,.&
    \end{array}$$
But since $v_{0}$ is killed by the action of the root-vectors
corresponding to positive roots in $\Phi_{G}^{+}$, one gets $$
  \begin{array}{rll}
    E_{i}'\cdot(g_{0}\cdot v_{0}) &
    =\frac{d}{dt}\Big(\Big(g_{0}g_{0}^{-1}\exp(t\,E_{i}') g_{0}\Big)\cdot v_{0}\Big)_{|t=0}
    \\&=\frac{d}{dt}\Big(\Big(g_{0}\exp\Big(t\,\mathrm{Ad}(g_{0}^{-1})\cdot E_{i}'\Big) \Big)\cdot
    v_{0}\Big)_{|t=0}\\&= g_{0}\cdot\Big( E_{i}\cdot  v_{0}\Big)\\
    &=0\,.
\end{array}$$
Hence the result.
\end{proof}
From the result~(\ref{spec}), we may then conclude:
\begin{lem}
$$2\,\|\beta_{w_{0}}\|^{2}+n/8\,,$$ is an eigenvalue of the square
of the Dirac operator.
\end{lem}
\begin{proof}
By the Freudenthal's formula, the Casimir eigenvalue
$c_{\gamma_{0}}$ of the representation $(\rho_{0},V_{0})$ is given
by
$$\|\beta_{w_{0}}^{G}+\delta_{G}\|^{2}-\|\delta_{G}\|^{2}=3\,\|\delta_{G}\|^{2}
+\|\delta_{K}\|^{2}-4\,<w_{0}\cdot\delta_{G},\delta_{K}>\,.$$ On
the other hand $$\|\beta_{w_{0}}\|^{2}=\|\delta_{G}\|^{2}
+\|\delta_{K}\|^{2}-2\,<w_{0}\cdot\delta_{G},\delta_{K}>\,.$$
Hence
$$c_{\gamma_{0}}=2\,\|\beta_{w_{0}}\|^{2}+\|\delta_{G}\|^{2}-\|\delta_{K}\|^{2}\,.$$
Now, the Casimir operator of $\gk$ acts on the spin representation
$\rho_{K}$ as scalar multiplication by $
\|\delta_{G}\|^{2}-\|\delta_{K}\|^{2}$, (cf.~lemma~2.2 in
\cite{Par}). Indeed, each dominant weight of $\rho_{K}$ being of
the form $w\cdot\delta_{G}-\delta_{K}$, $w\in W$, the eigenvalue
of the Casimir operator on each irreducible component is given by:
$$\|(w\cdot\delta_{G}-\delta_{K})+\delta_{K}\|^{2}-\|\delta_{K}\|^{2}=
\|w\cdot\delta_{G}\|^{2}-\|\delta_{K}\|^{2}=\|\delta_{G}\|^{2}-\|\delta_{K}\|^{2}\,.$$
On the other hand, the proof of the formula (\ref{Par}) shows that
the Casimir operator of $\gk$ acts on the spin representation
$\rho_{K}$ as scalar multiplication by
$\frac{\mathrm{Scal}}{8}=n/16$ (cf. \cite{Sul}), hence
\begin{equation}\label{sf}
\|\delta_{G}\|^{2}-\|\delta_{K}\|^{2}=n/16\,.
\end{equation} So
$$c_{\gamma_{0}}+n/16=2\,\|\beta_{w_{0}}\|^{2}+n/8\,.$$
\end{proof}
In order to conclude, we have to prove that
\begin{lem}
$$2\,\|\beta_{w_{0}}\|^{2}+n/8\,,$$ is the lowest eigenvalue of
the square of the Dirac operator.
\end{lem}
\begin{proof} Let $\gamma\in\widehat{G}$ such that there exists
$w\in W$ such that
$\mathrm{m}({\rho_{\gamma}}_{|K},\rho_{K,w})\geq 1$. Let
$\beta_{\gamma}$ be the dominant weight of ${\rho_{\gamma}}$.
First, since the Weyl group permutes the weights of
${\rho_{\gamma}}$,
$w^{-1}\cdot\beta_{w}=\delta_{G}-w^{-1}\cdot\delta_{K}$ is a
weight of ${\rho_{\gamma}}$. Hence
$$\|\beta_{\gamma}+\delta_{G}\|^{2}\geq
\|w^{-1}\cdot\beta_{w}+\delta_{G}\|^{2}\,,$$ (cf. for instance
Lemma~C, \S 13.4 in \cite{Hum}). So, from the Freudenthal formula,
$$c_{\gamma}=\|\beta_{\gamma}+\delta_{G}\|^{2}-\|\delta_{G}\|^{2}\geq
\|w^{-1}\cdot\beta_{w}+\delta_{G}\|^{2}-\|\delta_{G}\|^{2}\,.$$
But, using (\ref{sf})
$$\|w^{-1}\cdot\beta_{w}+\delta_{G}\|^{2}-\|\delta_{G}\|^{2}= 2\,
\|\beta_{w}\|^{2}+\|\delta_{G}\|^{2}-\|\delta_{K}\|^{2}=2\,
\|\beta_{w}\|^2+n/16\,.$$ Hence by the definition of
$\beta_{w_{0}}$, $$c_{\gamma}\geq 2\, \|\beta_{w}\|^2+n/16\geq
2\,\|\beta_{w_{0}}\|^2+n/16\,.$$ Hence the result.
\end{proof}

\end{document}